\documentclass{amsart}

\newtheorem{theorem}{Theorem}[section]
\newtheorem{lemma}[theorem]{Lemma}
\newtheorem{proposition}[theorem]{Proposition}
\newtheorem{corollary}[theorem]{Corollary}

\theoremstyle{definition}

\newtheorem{remark}[theorem]{Remark}
\newtheorem*{acknowledgement}{Acknowledgement}

\theoremstyle{remark}

\newcommand\mylabel[1]{\label{#1}}

\newcommand{\ZZ}{\mathbb{Z}}
\newcommand{\QQ}{\mathbb{Q}}

\newcommand{\FF}{\mathbb{F}}
\newcommand{\PP}{\mathbb{P}}
\renewcommand{\AA}{\mathbb{A}}
\newcommand{\GG}{\mathbb{G}}
\newcommand{\HH}{\mathbb{H}}

\newcommand  {\shN}     {\mathcal{N}}
\newcommand  {\shL}     {\mathcal{L}}

\newcommand  {\shTor}   {\mathcal{T}\!\text{\textit{or}}\,}

\newcommand  {\foY}     {\mathfrak{Y}}
\newcommand  {\foZ}     {\mathfrak{Z}}

\newcommand  {\Aut}     {\operatorname{Aut}}

\newcommand  {\Br}      {\operatorname{Br}}

\newcommand  {\disc}    {\operatorname{disc}}

\newcommand  {\Gal}     {\operatorname{Gal}}

\newcommand  {\gr}      {\operatorname{gr}}

\newcommand  {\Hom}     {\operatorname{Hom}}

\newcommand  {\hgt}     {\operatorname{ht}}

\newcommand  {\id}      {\operatorname{id}}

\renewcommand  {\k}     {\kappa}

\newcommand  {\Lie}     {\operatorname{Lie}}

\newcommand  {\lra}     {\longrightarrow}

\newcommand  {\maxid}   {\mathfrak{m}}

\newcommand  {\NS}      {\operatorname{NS}}

\renewcommand{\O}       {\mathcal{O}}

\newcommand  {\ord}     {\operatorname{ord}}

\newcommand  {\Pic}     {\operatorname{Pic}}

\newcommand  {\quadand} {\quad\text{and}\quad}

\newcommand  {\ra}      {\rightarrow}

\newcommand  {\Spec}    {\operatorname{Spec}}
\newcommand  {\Spf}     {\operatorname{Spf}}

\def\mydate{\number\day\space\ifcase\month \or January\or February\or March\or 
April\or May\or June\or July\or
August\or September\or October\or November\or December\fi \space\number\year}

\usepackage{amscd,amssymb}

\begin{document}

\title[Calabi--Yau threefolds]
      {Some Calabi--Yau threefolds with obstructed deformations
over the Witt vectors}

\author[Stefan Schroer]{Stefan Schr\"oer}
\address{Mathematische Fakult\"at, Ruhr-Universit\"at, 
         44780 Bochum, Germany}
\curraddr{Mathematisches Institut, Universit\"at K\"oln, Weyertal 86-90,
          50931 K\"oln, Germany}
\email{s.schroeer@ruhr-uni-bochum.de}

\subjclass{14J28, 14J30, 14J32,14K10}

\dedicatory{Revised version, 24 June 2003}

\begin{abstract}
I construct some
smooth Calabi--Yau threefolds in characteristic two and three
that do not lift to  characteristic zero. These threefolds
are pencils of supersingular K3-surfaces. The construction depends on
Moret-Bailly's pencil of
abelian surfaces and Katsura's analysis of generalized Kummer surfaces.
The threefold in characteristic two turns out to be nonrigid.
\end{abstract}

\maketitle
\tableofcontents

\section*{Introduction}

A fundamental result in deformation theory states that
each infinitesimal deformation of a compact complex K\"ahler
manifold $Y$ with trivial dualizing sheaf extends to a deformation
of arbitrarily high order. In other words, the base of the
versal deformation is smooth. This result is due
to several authors (see Bogomolov \cite{Bogomolov 1978}, 
Kawamata \cite{Kawamata 1992}, Ran \cite{Ran 1992},
Tian \cite{Tian 1987}, Todorov \cite{Todorov 1980}).

It is open whether or to what extend this holds true over ground fields
$k$ of characteristic $p>0$. I gave a sufficient condition
for unobstructedness involving crystalline cohomology groups
and formal liftings to characteristic zero \cite{Schroeer 2003}.
Deligne \cite{Deligne 1976} proved that any K3-surface lifts projectively to 
characteristic zero. In contrast,
Hirokado \cite{Hirokado 1999} found
a smooth threefold in characteristic $p=3$ that does not admit
a formal lifting. Here  the deformations over the
Witt vectors $W=W(k)$ are necessarily obstructed.

A natural question to ask: What is the relation between 
obstructions for deformations over Artin $k$-algebras and deformations
over Artin $W$-algebras? 
Unfortunately,  there are almost no examples 
that could shed some light on the situation.
The goal of this paper is construct some interesting new examples
of Calabi--Yau threefolds $Y$ in characteristic $p=2$ and $p=3$ that do not 
formally lift to characteristic zero.

Our examples are 
pencils of supersingular K3-surfaces over the projective line.
The construction draws from the theory of supersingular abelian surfaces
and generalized Kummer surfaces.
The starting point is Moret-Bailly's nontrivial pencil of
abelian surface \cite{Moret-Bailly 1981}.
It turns out that my threefold $Y$ in characteristic $p=2$ has some
nontrivial deformations over $k$-algebras.
This seems to be the first example of nonliftable Calabi--Yau
manifolds with nontrivial deformations.
It would be interesting to know whether or not
the Hodge--de Rham spectral sequence 
$H^r(Y,\Omega_Y^s)\Rightarrow \HH^{r+s}(Y,\Omega^\bullet_Y)$ degenerates.

The paper is organized as follows:
In the first section I consider threefolds with trivial canonical
class that are fibered in supersingular K3-surfaces. 
At this stage we do not care about  existence and merely deduce
some properties from our  assumptions.
In Section 2 we shall see that it is impossible to lift such threefolds
to characteristic zero, even in a formal way.
We turn in Section 3 to the existence question.
Here we use an approach of Moret-Bailly to construct threefolds
with trivial canonical class that are fibered in supersingular
abelian surfaces.
The construction works precisely at the primes $p=2$ and $p=3$,
but the two cases have different features.
This becomes clear in Section 4: Here the associated family of Kummer surfaces
easily gives examples, in characteristic $p=3$, 
of our nonliftable Calabi--Yau 
threefolds. This approach fails, however, in characteristic $p=2$
due to results of Shioda and Katsura.
We deal with this problem in Section 5:
Discarding Kummer surfaces, we work with generalized Kummer surfaces
by replacing the sign involution
by an automorphism of order three.
We compute the discriminant of the intersection form in Section 6,
which relies on descent theory and theta groups.
It turns out that  discriminant  jumps in pencils.
This is crucial in Section 7: Here we construct the desired nonliftable
Calabi--Yau threefold for $p=2$.
The nontrivial deformations are detected by jumping discriminants.
In the last Section we compute the Artin--Mazur formal group $\Phi^3_Y$.

\begin{acknowledgement}
I wish to thank Stefan M\"uller-Stach, Klaus K\"unneman,
Uwe Jannsen, and Jean-Marc Fontaine for helpful discussions.
I am also grateful to the referee,
who suggested several simplifications and 
pointed out some errors in the
first version.

\end{acknowledgement}

\section{Calculation of invariants}
\mylabel{calculation}

Let $k$ be an algebraically closed ground field of characteristic
$p>0$. Suppose $Y$ is a smooth proper threefold with $\omega_Y=\O_Y$, 
endowed with
a smooth morphism $f:Y\ra\PP^1$ such that all geometric
fibers $Y_{\bar{t}}$, $t\in\PP^1$ are 
K3-surfaces with Picard number $\rho(Y_{\bar{t}})=22$.
We shall see in the following sections that such threefold
actually exists in characteristics $p=2$ and $p=3$.
In this section I merely collect some  consequences from
the assumptions. 

We start by computing the Hodge numbers $h^{0i}(Y)=h^{i}(\O_Y)$:

\begin{proposition}
\mylabel{h invariants}
We have $h^{0}(\O_Y)=h^{3}(\O_Y)=1$ and $h^{1}(\O_Y)=h^{2}(\O_Y)=0$.
\end{proposition}

\proof
Clearly $h^0(\O_Y)=1$ because $Y$ is connected.
The Leray--Serre spectral sequence for $f:Y\ra\PP^1$ gives an exact sequence
$$
H^1(\PP^1,f_*\O_Y)\lra H^1(Y,\O_Y)\lra H^0(\PP^1,R^1f_*\O_Y).
$$
The fibers  $Y_t$ are connected, so $\O_{\PP^1}\ra f_*\O_Y$ is bijective, 
and the term on the left vanishes.
The cohomology groups $H^1(Y_t,\O_{Y_t})$ vanish for all $t\in\PP^1$, 
hence $R^1f_*\O_Y=0$ by semicontinuity,
so the term on the right vanishes as well. We conclude $h^1(\O_Y)=0$.
The remaining Hodge numbers follow from Serre duality.
\qed

\medskip
Smooth proper $n$-folds $X$ with $\omega_X=\O_X$ 
and $h^1(\O_X)=\ldots=h^{n-1}(\O_X)=0$ are called 
\emph{Calabi--Yau} $n$-folds.
They play an important role in the
classification of schemes with trivial canonical class
\cite{Beauville 1983} and mirror symmetry. 
Our $Y$ is an example of a Calabi--Yau threefold.

\begin{proposition}
\mylabel{fundamental group}
The algebraic fundamental group $\pi_1(Y,y)$ vanishes.
\end{proposition}

\proof
Let $y\in Y$ be a closed point and $z\in\PP^1$ be its image.
According to \cite{SGA 1}, Expos\'e X,
Corollary 1.4, we have an exact sequence
$$
\pi_1(Y_z,y)\lra \pi_1(Y,y)\lra\pi_1(\PP^1,z).
$$
The K3-surface $Y_z$ and the projective line $\PP^1$ are simply connected,
so the same holds for $Y$.
\qed

\medskip
Let $\eta\in\PP^1$ be the generic point. Then $Y_\eta$ is a K3-surface
over the function field $\k(\eta)$. Choose an algebraic closure
of $\k(\eta)$ and let $Y_{\bar{\eta}}$ 
be the corresponding geometric generic fiber.

\begin{lemma}
\mylabel{generic picard}
The Picard group $\Pic(Y_\eta)$ is free of rank $\rho(Y_\eta)=22$, and the
inclusion map $\Pic(Y_\eta)\ra\Pic(Y_{\bar{\eta}})$ is bijective.
\end{lemma}

\proof
To simplify notation, set $S=Y_\eta$ and $K=\k(\eta)$, 
let $K\subset \tilde{K}\subset\bar{K}$ be separable and algebraic
closures, and write $\tilde{S}=S\otimes\tilde{K}$ and 
$\bar{S}=S\otimes\bar{K}$.
First note that the map $\Pic(\tilde{S})\ra\Pic(\bar{S})$ 
is bijective, because the Picard scheme
$\Pic_{S/K}$ is smooth and 0-dimensional.

According to Tsen's Theorem, we have $\Br(K)=0$.
This implies that the injection 
$\Pic(S)\subset\Pic(S/K)$ is bijective, where
the right hand side denotes the group of rational points in 
$\Pic_{S/K}$.
Note that $\Pic(S/K)=\Pic(\tilde{S})^G$, where $G=\Gal(\tilde{K}/K)$ 
is the Galois group.
By  assumption,  $\Pic(\tilde{S})=\Pic(\bar{S})$ is free of rank 22,
so it remains to check that $G$ acts trivially on $\Pic(\tilde{S})$.
The following argument suggested by the referee proves this:

Fix a prime power $l^n$ with $l\neq p$ and consider the
constructible sheaf $\mu_{l^n}$ on $Y$, which is constant.
Its direct image $R^2 f_*(\mu_{l^n})$ is a constructible  sheaf,
and locally constant because $f:Y\ra \PP^1$ is smooth
(\cite {Freitag; Kiehl 1988}, Chapter I, Lemma 8.13).
Then it must be constant because $\PP^1$ is simply connected.
Passing to the generic fiber, we infer that  
the $G$-action on $H^2(\tilde{S},\mu_{l^n})$ is trivial.
Using the Kummer sequence, we deduce that  $G$ acts trivially on
$\Pic(\tilde{S})$  as well.
\qed

\begin{proposition}
\mylabel{picard group}
The Picard group $\Pic(Y)$ is a free group of rank $\rho(Y)=23$.
\end{proposition}

\proof
The closure of any Cartier divisor on $Y_\eta$ is a Cartier divisor
on $Y$, because $Y$ is locally factorial.
Therefore the restriction map $\Pic(Y)\ra\Pic(Y_\eta)$ is surjective. 
An invertible $\O_Y$-module
that is trivial on $Y_\eta$ is trivial on some open neighborhood
of $Y_\eta$. It is  then trivial on the preimage of some
nonempty open subset $U\subset\PP^1$, because
$f:Y\ra\PP^1$ is a closed map.  Hence it comes from a Cartier divisor
supported by the preimage of $\PP^1-U$. In other words, the sequence
$$
0\lra\Pic(\PP^1)\lra\Pic(Y)\lra\Pic(Y_\eta)\lra 0
$$ 
is exact. By Lemma \ref{generic picard}, $\Pic(Y_\eta)$ is a free group 
of rank $\rho(Y_\eta)=22$, and the
statement follows.
\qed

\medskip
We now turn to the \emph{Brauer group} $\Br(Y)=H^2(Y,\GG_m)$, where
the cohomology is taken with respect to the \'etale topology.

\begin{proposition}
\mylabel{brauer group}
The Brauer group $\Br(Y)$ is annihilated by $p$.
\end{proposition}

\proof
Set $S=Y_\eta$ and let $\bar{S}=Y_{\bar{\eta}}$ the geometric generic fiber.
It follows from \cite{GB II}, Corollary 1.8 that the restriction map
$\Br(Y)\ra\Br(S)$ is injective.
The Leray--Serre spectral sequence for the structure morphism
$g:S\ra\Spec\k(\eta)$ 
yields a spectral sequence 
$H^r(\eta,R^sg_*\GG_{m,S})\Rightarrow H^{r+s}(S,\GG_{m,S})$.
The stalk $(R^1g_*\GG_{m,S})_{\bar{\eta}}=\Pic(\bar{S})$ 
is a finitely generated free group of finite rank,
so $H^1(\eta,R^1g_*\GG_{m,S})=0$.
The group $H^2(\eta,g_*\GG_{m,S})=\Br(\eta)$ vanishes by Tsen's Theorem.
It follows that the edge map $\Br(S)\ra H^0(\eta,R^2g_*\GG_{m,S})$ 
is injective. The latter group is contained in $\Br(\bar{S})$.
Now $\bar{S}$ is a K3-surface with Picard number $\rho=22$ 
over an algebraically closed field.
Artin's argument for \cite{Artin 1974}, Theorem 1.7 shows that the abelian
group 
$\Br(\bar{S})$ is annihilated by $p$.
\qed

\medskip
Next we  calculate the $l$-adic Betti numbers
$b_i(Y)=\dim_{\QQ_l} H^i(Y,\QQ_l)$ and the Euler characteristic
$e(Y)=\sum(-1)^ib_i(Y)$, 
where $l$ is any prime number different from $p$.

\begin{proposition}
\mylabel{betti numbers}
The $l$-adic Betti numbers $b_i=b_i(Y)$ and Euler characteristic
are as follows:
$b_0=b_6=1$, and $b_1=b_3=b_5=0$, and $b_2=b_4=23$, and $e=48$.
\end{proposition}

\proof
Obviously  $b_0=1$ because $Y$ is connected. 
Proposition \ref{fundamental group}
implies that $H^1(Y,\ZZ/l^n\ZZ)=0$ for each $n>0$, hence $b_1=0$.
To proceed, consider the Leray--Serre spectral sequence
$$
H^r(\PP^1,R^sf_*(\mu_{l^n}))\Longrightarrow H^{r+s}(Y,\mu_{l^n}).
$$
The constructible sheaves $R^sf_*(\mu_{l^n})$ 
are constant, by the same argument
as in the proof of Lemma \ref{generic picard}.
It follows from the proper base change theorem 
that $R^sf_*(\mu_{l^n})=0$  for $s$ odd, because then the K3-surfaces
$Y_t$ satisfy
$H^s(Y_t,\mu_{l^n})=0$. Therefore the spectral sequence must degenerate.
Moreover, $H^r(\PP^1,R^sf_*(\mu_{l^n}))$ vanishes for $r$ odd
because $\PP^1$ is simply connected.

From this we immediately deduce $H^3(Y,\mu_{l^n})=0$
and hence $b_3=0$.
The module $H^2(Y,\mu_{l^n})$ 
is isomorphic to the sum of $H^2(\PP^1,\mu_{l^n})$, which has rank $1$,
and $H^0(\PP^1,R^2f_*(\mu_{l^n}))$, which has rank $22$.
The latter holds because by  the stalks  
$R^2f_*(\mu_{l^n})_t=H^2(Y_t,\mu_{l^n})$ have rank
$22$.
Consequently $b_2=23$, and the remaining Betti numbers follow from
Poincar\'e duality.
\qed

\medskip
We finally observe that $Y$ carries an ample invertible sheaf:

\begin{proposition}
\mylabel{projective}
The scheme $Y$ is projective.
\end{proposition}

\proof
It suffices to check that the morphism $f:Y\ra\PP^1$ is projective.
The generic fiber $Y_\eta$ is a regular surface, hence there is an ample
invertible $\O_{Y_\eta}$-module $\shL_\eta$. 
Extend it to an invertible $\O_Y$-module $\shL$ and fix a point $t\in\PP^1$.
The specialization map $\Pic(Y_\eta)\ra\Pic(Y_t)$ 
has finite cokernel, because all fibers
have Picard number $\rho=22$. Hence $\shL\cdot C>0$ for any curve 
$C\subset Y_t$.
We also have $\shL_t\cdot\shL_t=\shL_\eta\cdot\shL_\eta>0$.
This implies that $\shL_t=\shL|_{Y_t}$ is ample, and therefore $\shL$ 
is $f$-ample.
\qed

\section{Nonliftability}

We keep the assumptions from the preceding section.
We say that our Calabi--Yau threefold
$Y$ admits a \emph{formal lifting to characteristic zero}
if there is a complete local noetherian ring $R$ with residue
field $R/\maxid=k$ containing $\ZZ$ as a subring,
 together with a proper flat formal $R$-scheme $\foZ\ra\Spf(R)$
with closed fiber $\foZ_0=Y$.
For example, a smooth proper scheme $X$ 
admits a formal lifting to characteristic zero
if $H^2(X,\Theta_X)=0$. The first main result of this paper is:

\begin{theorem}
\mylabel{no lifting}
The Calabi--Yau threefold $Y$ does not admit a formal lifting
to characteristic zero.
\end{theorem}

\proof
Suppose to the contrary that there is a formal lifting $\foZ\ra\Spf(R)$ 
to characteristic zero. We may assume that
$R$ is integral, say with field of fractions $R\subset Q$.
Choose a descending sequence of ideals $\maxid=I_0\supset I_1\supset\ldots$
so that $I_n/I_{n+1}$ have length one.
Let $Y_n\subset \foZ$ be the infinitesimal neighborhood of $Y$
defined by the ideal $I_{n+1}\O_{\foZ}$.
The exact sequence
$$
H^1(Y,\O_Y)\ra\Pic(Y_{n+1})\lra\Pic(Y_n)\lra H^2(Y,\O_Y)
$$
implies that the restriction map $\Pic(\foZ)\ra\Pic(Y)$ is bijective.
According to Proposition \ref{projective}, there is an ample invertible
$\O_Y$-module,
so Grothendieck's Algebraization Theorem tells us that $\foZ$ is
algebraizable. Write $\foZ=Z^\wedge$ 
for some smooth proper $g:Z\ra \Spec(R)$.
Note that $\Pic(Z)=\Pic(\foZ)=\Pic(Y)$. 
In particular, the dualizing sheaf $\omega_{Z/R}$ is trivial,
and the generic fiber $Z_\eta$ over $R$ is a Calabi--Yau threefold
in characteristic zero.

Now Hirokado's argument \cite{Hirokado 1999} yields a contradiction:
Using that the constructible sheaves $R^3g_*(\mu_{l^n})$ are constant,
we deduce $b_3(Z_{\bar{\eta}})=b_3(Y)=0$.
On the other hand, the Hodge--de Rham spectral sequence
for $Z_{\bar{\eta}}$ degenerates, hence
$b_3(Z_{\bar{\eta}})\geq h^{3,0}(Z_{\bar{\eta}})=1$, contradiction.
\qed

\medskip
As a consequence, the deformation functor $D_Y$ is not formally smooth.
Here $D_Y$ is the functor  which associates to an
Artin algebra $R$ over the Witt vectors $W(k)$ 
the set $D_Y(R)$ of isomorphism classes of proper flat
$R$-schemes $\foY$ together with an identification $\foY_0\simeq Y$.

\begin{remark}
The 1-dimensional Calabi--Yau manifolds are  the elliptic curves,
and Grothendieck showed that any smooth proper curve lifts to characteristic
zero (\cite{SGA 1}, Expos\'e III, Corollary 7.4).
The 2-dimensional Calabi--Yau manifolds are the K3-surface.
According to Deligne \cite{Deligne 1976}, every K3-surface lifts
to characteristic zero.
Hirokado's \cite{Hirokado 1999} counterexample showed that one 
cannot hope for any general result
like this in higher dimensions.
The upshot of this paper is that nonliftability  might be quite frequent
for Calabi--Yau manifolds.
\end{remark}

\section{Moret-Bailly's pencil of abelian surfaces}
\mylabel{abelian pencil}

The task now is to construct pencils of K3-surfaces $Y$
as in Section \ref{calculation}.
In this section we take a first step into that direction by
constructing a pencil of abelian surfaces with similar properties.

Fix an algebraically closed field  $k$  of characteristic $p>0$.

Let $A$ be the \emph{superspecial} abelian surface, that is, $A$ is
a 2-dimensional abelian variety isomorphic to the product of two
supersingular elliptic curves. Note that $A$ does not depend, up to
isomorphism, on the choice of the supersingular elliptic curves
by \cite{Shioda 1979}, Theorem 3.5.
Its \emph{$a$-number}
$$
a(A)=\dim_k\Hom_{\gr}(\alpha_p,A)
$$
equals $a(A)=2$, and we have an embedding of group schemes 
$\alpha_p^{\oplus 2}\subset A$. Actually, Oort proved in
\cite{Oort 1975}, Theorem 2 that an abelian surface
is superspecial if and only if its $a$-number is $a=2$.

Fix an integer $n\geq 1$ and an exact sequence
\begin{equation}
\label{extension}
0\lra \O_{\PP^1}(-n)\lra\O_{\PP^1}^{\oplus 2} \lra \O_{\PP^1}(n)\lra 0.
\end{equation}
Such extensions correspond to  pairs of sections 
$r,s\in H^0(\PP^1,\O_{\PP^1}(n))$ without  common zeros.
Let $X'=A\times\PP^1$ be 
the constant relative abelian variety over $\PP^1$,
and $H\subset X'$ the relative radical subgroup scheme of height $\hgt(H)=1$ 
whose $p$-Lie algebra
is isomorphic to 
$$
\Lie(H/\PP^1)=\O_{\PP^1}(-n)\subset
\O_{\PP^1}^{\oplus 2} = \Lie(\alpha_p^{\oplus 2}/\PP^1)=\Lie(X'/\PP^1),
$$
as explained in \cite{SGA 3a}, Expos\'e VII$_A$,  Theorem 7.2.
Here the inclusion 
$\O_{\PP^1}(-n)\subset\O_{\PP^1}^{\oplus 2} $ comes from 
the extension (\ref{extension}). This works because both
  Lie bracket and   $p$-th power operation are trivial for
$\Lie(\alpha_p^{\oplus 2})$.

The quotient $X=X'/H$ is a relative abelian surface over $\PP^1$, whose
fibers $X_t$, $t\in\PP^1$ are supersingular abelian surfaces.
This construction is due to Moret-Bailly \cite{Moret-Bailly 1979}, who
considered the case $n=1$.

\begin{proposition}
\mylabel{lie algebra}
We have $\Lie(X/\PP^1)\simeq\O_{\PP^1}(-np)\oplus\O_{\PP^1}(n)$.
\end{proposition}

\proof
Moret-Bailly showed this for $n=1$. In the general case, the family
$X\ra\PP^1$ comes from the family with $n=1$  via pullback along
the map $(r,s):\PP^1\ra\PP^1$ of degree $n$.
For the convenience of the reader I recall Moret-Bailly's arguments:

First note that $\Lie(X/\PP^1)$ is dual to the conormal sheaf 
$N_{H/X'}$ restricted to the zero section $\PP^1\subset H$.
To see this we use cotangent complexes (see
\cite{SGA 6}, Expos\'e VIII, and \cite{Illusie 1971}).
Consider the commutative diagram of cotangent complexes
$$
\begin{CD}
L^\bullet_{X/\PP^1}|_H @>>> L^\bullet_{H/\PP^1} @>>> L^\bullet_{H/X}\\
@. @A\simeq AA @AA\id A \\
@. L^\bullet_{X'/X}|_H@>>> L^\bullet_{H/X}@>>> L^\bullet_{H/X'}
\end{CD}
$$
in the derived category  of perfect complexes on $H$.
Here the rows are distinguished triangles coming from the two
compositions $H\ra X\ra\PP^1$ and $H\ra X'\ra X$.
The vertical map on the left is the base change map induced
from the cartesian diagram
$$
\begin{CD}
H @>>> X'\\
@VVV @VVV\\
\PP^1 @>> 0> X
\end{CD}
$$
This base change map is a quasiisomorphism 
according to \cite{Illusie 1971}, Chapter II, Corollary 2.2.3,
because $\shTor^{\O_X}_r(\O_{X'},\O_{\PP^1})=0$ for $r>0$ since 
$X'\ra X$ is flat.
It then follows that there is a quasiisomorphism
$L^\bullet_{X/\PP^1}|_H\ra L^\bullet_{H/X'}[-1]$ in the derived category, 
and the induced bijection on
cohomology yields the desired identification
$\Omega^1_{X/\PP^1}|_H\simeq N_{H/X'}$.

To calculate the latter, we may replace $X'$ by $X''=\GG_{a,\PP^1}^{\oplus 2}$
via the canonical composition of closed embeddings
$$
\alpha_{p,\PP^1}^{\oplus 2}\subset {}_pX'\subset \GG_{a,\PP^1}^{\oplus 2}
$$
as explained in \cite{Moret-Bailly 1979}, page 136.
Let $U,V\in H^0(\PP^1,\O_{\PP^1}(1))$ be homogeneous coordinates, 
and write $\GG_a^{\oplus 2}=\Spec k[x,y]$.
Let $R=k[x,y,U,V]$ be the homogeneous coordinate ring for $X''$, and 
$I\subset R$ the
homogeneous ideal for $H\subset X''$. Then
$I=(y^p,x^p,sx-ry)$, where $r,s$ are the two homogeneous polynomials 
in $U,V$ of degree $n$
defining $\Lie(H)\subset\Lie(X'/\PP^1)=\Lie(X''/\PP^1)$.
 The  Hilbert--Burch Theorem
(see \cite{Eisenbud 1995}, Theorem 20.12) tells us that 
\begin{equation}
\label{hilbert-burch}
0\lra R(-n)\oplus R(-np)\stackrel{\varphi_2}{\lra}
R\oplus R\oplus R(-n)\stackrel{\varphi_1}{\lra} I\lra 0
\end{equation}
is an exact sequence of graded $R$-modules  
where the maps are given by the matrices
$$
\varphi_2=
\left(\begin{array}{cc}
0 & r^p\\
ry-sx & -s^p\\
x^p & (ry-sx)^{p-1}
\end{array}\right) 
\quadand
\varphi_1=
\left(\begin{array}{ccc}
y^p, & x^p, & sx-ry
\end{array}\right).
$$ 
Indeed, the entries of $-r^p\varphi_1$ are the $2\times 2$-minors of 
$\varphi_2$ with the
appropriate signs
(note that the  sequence 
in \cite{Moret-Bailly 1979} on page 137 corresponding to (\ref{hilbert-burch})
is slightly incorrect).
Tensoring the resolution with the residue ring
$\bar{R}=R/(x,y)$, we obtain an exact sequence
$$
\bar{R}(-n)\oplus \bar{R}(-np)\stackrel{\bar{\varphi}_2}{\lra}
\bar{R}\oplus \bar{R}\oplus \bar{R}(-n)\stackrel{\bar{\varphi}_1}{\lra} 
I/(x,y)I\lra 0,
$$
which in turn gives an exact sequence
$$
0\lra \bar{R}(-np)\lra
\bar{R}\oplus \bar{R}\oplus \bar{R}(-n)\stackrel{\bar{\varphi}_1}{\lra} 
I/(x,y)I\lra 0.
$$
The map on the left is defined by the transpose
of $(r^p,-s^p,0)$, and its cokernel is the graded $\bar{R}$-module  
$\bar{R}(np)\oplus\bar{R}(-n)$.
The letter defines the locally free $\O_{\PP^1}$-module
$\O_{\PP^1}(np)\oplus\O_{\PP^1}(-n)$, which by construction 
is dual to $\Lie(X/\PP^1)$.
\qed

\medskip
This has an interesting consequence for the dualizing sheaf $\omega_{X}$:

\begin{corollary}
\mylabel{dualizing sheaf}
The dualizing sheaf $\omega_{X}$ is trivial
if and only if either $p=3$ and $n=1$, or $p=2$ and $n=2$.
\end{corollary}

\proof
We have $\omega_{X/\PP^1}=\O_X(np-n)$ by Proposition \ref{lie algebra} and
$\omega_{\PP^1}=\O_{\PP^1}(-2)$, so $\omega_X=\O_X(np-n-2)$.
Consequently $\omega_X$ is trivial if and only if $n(p-1)=2$,
and the result follows.
\qed

\medskip
We shall see in the following sections how to obtain from
$X$  as above some interesting Calabi--Yau threefolds.

\section{A pencil of Kummer surfaces}
\mylabel{kummer pencil}

We now examine the family $X\ra\PP^1$ of abelian surface constructed
in the preceding section in the special case $p=3$ and $n=1$.
This is the first of the two cases where $\omega_X$ is trivial.

Since $p\neq 2$ and $k$ is algebraically closed, the kernel ${}_2 X\subset X$
of the multiplication by two map $[2]:X\ra X$ consists of sixteen disjoint
sections. We may also view ${}_2X$ as the relative fixed point scheme
of the sign change map $[-1]:X\ra X$. Let $Z=X/[-1]$ be 
the corresponding quotient
scheme. It inherits a fibration $Z\ra\PP^1$ whose fibers $Z_t$ are singular 
Kummer surfaces.
Let $D_\eta\subset Z_\eta$ be the reduced singular locus and $D\subset Z$ be
its closure. By flatness, each fiber $D_t$ is the reduced singular locus. 
The sixteen singular points on $Z_t$ are rational double points of type $A_1$.
The blowing up of $D_t$ yields the minimal resolution.
Let $X\ra Z$ be the blowing up with center $D\subset Z$. Then $Y$ is a smooth
threefold, and the  fibers $Y_t$ of the induced
family $Y\ra\PP^1$ are Kummer K3-surfaces.

\begin{proposition}
\mylabel{supersingular fibers}
We have  $\rho(Y_{\bar{t}})=22$ for every $t\in\PP^1$.
\end{proposition}

\proof
Since $E$ is supersingular, the Picard number of $A$ is $\rho(A)=6$, so 
$\rho(B_{\bar{t}})=6$ as well. It then follows from 
\cite{Shioda 1979}, Proposition 3.1 that $\rho(Y_{\bar{t}})=16+6=22$.
\qed

\begin{proposition}
\mylabel{omega trivial}
The dualizing sheaf $\omega_Y$ is trivial.
\end{proposition}

\proof
Let $\eta\in\PP^1$ be the generic point. Then $\omega_{Y_\eta}$ is trivial, 
because $Y_{\eta}$
is a K3-surface. Hence $\omega_Y=\O_Y(m)$ for some integer $m$.
The relative dualizing sheaf $\omega_{Y/Z}$ 
is supported by the exceptional locus
of $Y\ra Z$. It must be trivial, because it is trivial on the generic fiber.
We conclude $\omega_Z=\O_Z(m)$. 
The canonical map $f:X\ra Z$ is \'etale on  $X-{}_2X$, so
the invertible sheaves $\omega_X$ and $f^*(\omega_Z)$ coincide outside ${}_2X$.
Since $X$ satisfies Serre's condition $(S_2)$ and
${}_2X\subset X$ has codimension two, we must have 
$\omega_X\simeq f^*(\omega_Z)$  by
\cite{Hartshorne 1994}, Theorem 1.12. Consequently $m=0$, 
and $\omega_Y$ must be trivial.
\qed

\medskip
We conclude that $Y$ satisfies the assumptions from
Section \ref{calculation}, so $Y$ is a Calabi--Yau threefold
with $b_2=23$ that does not admit a formal lifting to characteristic zero.
 
Note that the preceding approach  does not
work for $p=2$, because then $Z_t$ is a  rational surfaces
with an elliptic singularity
(see \cite{Katsura 1987}, Corollary 2.12, or
\cite{Katsura 1977}, Theorem A, or \cite{Shioda 1974}, Proposition 1).
We shall sidestep this problem in the next  sections.

\section{Generalized Kummer surfaces}
\mylabel{generalized kummer}

In this section we construct some 
\emph{generalized Kummer surfaces}.
This terminology is due to Katsura \cite{Katsura 1987},
and denotes K3-surfaces arising as minimal resolutions 
of quotients of  abelian surfaces by  finite group actions.
Katsura studied mainly the
case $p\geq 3$.
Throughout  we work in characteristic $p=2$. Let $B$ be an abelian surface.
According to \cite{Katsura 1977}, Theorem A, the corresponding
singular Kummer surface $B/[-1]$ is birational to a K3-surface if and only if
$B$ is ordinary, and a rational surface  if and only if 
$B$ is supersingular. To obtain K3-surfaces from supersingular
abelian surfaces, we shall replace the sign involution $[-1]$ by
an automorphism of order three.

In this section our ground field $k$ contains $\FF_4$, but is
not necessarily algebraically closed.
Let $E$ be the supersingular elliptic curve given
by the Weierstrass equation $x^3=y^2+y$ and fix
a primitive third root of unity $\zeta\in\FF_4$.
Let $\varphi:E\ra E$ be the automorphism given by $(x,y)\mapsto(\zeta x, y)$.
This defines an action of $G=\ZZ/3\ZZ$ on $E$.
The induced representation on $H^0(E,\Omega^1_E)$ 
is scalar multiplication by $\zeta$,
because $dx$ is an invariant differential form on $E$.
The action on $E$ has three fixed points, namely $(0,0)$, $(0,1)$, 
and $\infty$.
Each fixed point $x\in E$ comes along with a  representations
$\rho_x:G\ra\Aut_{\k(x)}(\maxid_x/\maxid_x^2)$ 
taking values in $\FF_4^\times\subset\k(x)^\times$
The generator $\varphi$ maps  to   $\zeta\in\FF_4^\times$ under each $\rho_x$ 
because $\maxid_x/\maxid_x^2=\Omega_E^1(x)$.

Now consider the superspecial abelian surface $A=E\times E$, endowed with the
action of $G$ via $\phi=(\varphi,\varphi)$. Let $\alpha_2\subset A$ 
be any embedding of group schemes.

\begin{proposition}
\mylabel{commuting action}
The subgroup scheme $\alpha_2\subset A$ is stable under the 
$G$-action.
\end{proposition}

\proof
We have to check that $\phi(\alpha_2)=\alpha_2$ as subschemes in $A$.
The origin  $x\in A$ is clearly a fixed point for $\phi$.
Let $\maxid\subset\O_{A,x}$ be the maximal ideal.
The first order infinitesimal neighborhood of $x$ is the spectrum of
$R=k\oplus \maxid/\maxid^2$. 
As remarked above, $\phi$ acts via scalar multiplication by $\zeta$
on $\maxid/\maxid^2$. 
Hence $G$ leaves any subscheme of $\Spec(R)$ invariant.
\qed

\medskip
We see that  the $G$-action on the superspecial abelian surface
$A=E\times E$ descends to a $G$-action
on the supersingular abelian surface $B=A/\alpha_2$.
Its fixed points are easy to determine:
We saw that $\varphi$ has precisely three
fixed points. In turn, $\phi:A\ra A$ has 
$9=3\times 3$ fixed points, and the same holds on $B$.
These fixed points correspond to the singularities on $B/G$,
which is a proper normal surface.

\begin{proposition}
\mylabel{singularity type}
The singularities on $B/G$ are rational double points of type $A_2$.
\end{proposition}

\proof
It suffices to check this for the origin $x\in A$, because
we may replace any fixed point of $\varphi:E\ra E$ by the origin of $E$.
Let $\hat{E}$ be the formal completion at the origin $y\in E$, 
which is a formal group 
isomorphic to $\hat{\GG}_a$ because $E$ is supersingular.
Hence there uniformizer $u\in\O_{E,y}^\wedge$ such that the multiplication 
$\hat{E}\times\hat{E}\ra \hat{E}$ is
given by the formal group law $u\mapsto u+u'$.
In turn, there is a regular system of parameters $u,v\in\O_{A,x}^\wedge$  
such that $u\mapsto u+u'$ and $v\mapsto v+v'$ is the
formal group law for $\hat{A}$.

Write $\alpha_2=\Spec k[\epsilon]$.
The action $\alpha_2\times\hat{A}\ra\hat{A}$ 
is given by $u\mapsto u+i\epsilon$ and $v\mapsto v+j\epsilon$ 
for some $i,j\in k$ that do not vanish both, say $j\neq 0$.
A power series $f(u,v)$ is invariant under the $\alpha_2$-action if and only if
$f(u,v)=f(u+i\epsilon, v+j\epsilon)$.
For example, $x'=u^2$ and $y'=u+vi/j$ are $\alpha_2$-invariants.
Using that 
$\O_{A,x}^\wedge$ is a free module of rank two over $\O_{B,x}^\wedge$
by \cite {Serre 1965}, page IV-37, Proposition 22, 
we infer that $x',y'\in\O_{B,x}^\wedge$ is a regular system of
parameters.
The induced $G$-action has $\phi^*(x')\equiv\zeta^2 x'$ and 
$\phi^*(y')\equiv\zeta y'$ 
modulo $\maxid^2$, where $\maxid\subset \O_{B,x}^\wedge$
denotes the maximal ideal.
Using Lemma \ref{diagonizable} below, we find another
regular system of parameters $x,y\in\O_{B,x}^\wedge$ with 
$\phi^*(x)=\zeta^2 x$ and $\phi^*(y)=\zeta y$.
The $G$-invariants are then generated by  $x^3,xy,y^3$,  and this means that
$B/G$ acquires a rational double point of type $A_2$.
\qed

\medskip
A local computation shows that blowing-up the reduced singular locus
in a rational double point of type $A_2$ yields the minimal resolution,
and the exceptional divisor consists of two smooth rational $(-2)$-curves
intersecting transversely.

\begin{theorem}
\mylabel{k3 surface}
The minimal resolution $S$ of $B/G$ is a K3-surface.
\end{theorem}

\proof
The $G$-action on $B$ has the following four properties:
First, it has no fixed curves. Second, the generator of $G$ has
some isolated fixed points. Third, all singularities on $B/G$
are rational singularities by Proposition \ref{singularity type}.
Finally, the induced $G$-action on $H^0(B,\omega_B)=k$ is trivial.
To see this choose a fixed point $x\in B$. 
Let $\maxid\subset\O_{B,x}$ be the maximal ideal. We saw that
the $k$-linear $G$-action on $\maxid/\maxid^2$ is diagonizable
with eigenvalues $\zeta,\zeta^2$. Using $\det(\maxid/\maxid^2)=\Omega^2_B(x)$, 
we infer that $G$ must act via $1=\zeta\cdot\zeta^2$ on global 2-forms.

Assuming $p\neq 2$, Katsura showed that under the preceding four conditions
the minimal resolution $S$ of $B/G$
is a K3-surface (\cite{Katsura 1987}, Theorem 2.4 and Lemma 2.7).
Actually, he showed without any assumption on the characteristic
that $\omega_S$ is trivial, and that $S$ is
neither abelian, hyperelliptic, nor quasihyperelliptic.
It then follows from the classification of surfaces that
$S$ is either a K3-surface or a nonclassical Enriques surface
\cite{Bombieri; Mumford 1977}. The latter
exists only in characteristic $p=2$.
An Enriques surface, however, has by definition $l$-adic Betti number $b_2=10$.
On the other hand we have $b_2(S)> 18$, due to the exceptional
curves coming from the nine rational double points of type $A_2$ on $B/G$.
Therefore, $S$ must be a K3-surface. 
\qed

\medskip
We used the following fact in the proof for Proposition 
\ref{singularity type}.
I include a proof for lack of reference.

\begin{lemma}
\mylabel{diagonizable}
Let $R$ be a regular local noetherian $k$-algebra and
$f:R\ra R$ an $k$-automorphism of finite order $m>0$ prime to the 
characteristic $p$. It $k$ contains a primitive $m$-th root of unity $\zeta$,
then there is a regular system of parameters $u_1,\ldots,u_d\in R^\wedge$
and integers $n_1,\ldots,n_d$
such that  $f(u_i)=\zeta^{n_i} u_i$.
\end{lemma}

\proof
Let $\maxid\subset R$ be the maximal ideal. The induced $k$-linear map
on $\maxid/\maxid^2$ is diagonizable with eigenvalues of the form 
$\zeta^{n_i}$, 
because $m$ is prime to $p$ and $\zeta\in k$.
Hence there is a regular system of parameters $u_i\in R$ whose classes
modulo $\maxid^2$ are eigenvectors for the eigenvalues $\zeta^{n_i}$.
Suppose we have a homogeneous ideal $I\subset R$ with 
$f(u_i)\equiv \zeta^{n_i}u_i$ modulo $I$.
If $I\neq 0$, then
choose a homogeneous subideal $J\subsetneq I$ so that $I=J+ku^l$, 
with the usual multiindex notation
$u^l=u_1^{l_1}\ldots u_d^{l_d}$. Then
$f(u_i)\equiv \zeta^{n_i} + yu^l$
modulo $J$ for some $y\in k$. 
We distinguish two cases: If $m$ divides $e=\sum l_in_i$, then
the congruence
$$
u_i=f^m(u_i)\equiv\zeta^{mn_i}u_i + yu^l(1+\zeta^e+\ldots+\zeta^{e(m-1)})
=u_i + yu^l m
$$
modulo $J$ implies $y=0$. If $m$ does not divide $e$, then the congruence
$$
f(u_i+xu^l)\equiv \zeta^{n_i}u_i +yu^l +x\zeta^eu^l=
\zeta^{n_i}(u_i+xu^l) +u^l(y+x(\zeta^e-1))
$$
modulo $J$ shows that me may choose $x\in k$ so that $y+x(\zeta^e-1)=0$ 
and then replace
$u_i$ by $u_i+xu^l$. The principle of noetherian induction finishes the proof.
\qed

\section{Discriminants and Artin invariants}
\mylabel{discriminants}

We keep the notation as in the preceding section,
such that $A=E\times E$ and $B=A/\alpha_2$  for some embedding 
$\alpha_2\subset A$. We saw that the minimal resolution $S$ of $B/G$
is a K3-surface, so $\Pic(S)$ is a free group of rank $\rho(S)\leq 22$.
The goal of this section is to calculate the discriminant
of this intersection pairing on $\Pic(S)$. 
It turns out that it depends in a subtle
way on the position of $\alpha_2\subset A$, and this will be crucial in
Section \ref{nonrigidity}.

First, we examine the effect of the isogeny $g:A\ra B$ on Picard groups.
Consider the following four elliptic curves $C_i\subset A$:
Set $C_1=E\times 0$ and $C_2=0\times E$, 
let $C_3=\Delta_E$ be the diagonal subscheme,
and $C_4=\Gamma_\varphi$ the graph of the automorphisms $\varphi:E\ra E$.
Note that  $\Delta_E\cap\Gamma_\varphi$ is the fixed scheme
of $\varphi:E\ra E$, which consists of three reduced points, 
hence $C_3\cdot C_4=3$.

Let $V_A$ be the free abelian group generated by the $C_i$.
To simplify calculations, we introduce another basis $C_i'$ for $V_A$
as follows: Set $C_i'=C_i$ for $i=1,2$, and $C_i'=C_i-C_1-C_2$ for $i=3,4$.
The corresponding intersection matrix has an 
orthogonal decomposition
\begin{equation}
\label{intersection matrix}
(C_i'\cdot C_j')=
\left(
\begin{array}{cc}
0 & 1\\
1 & 0\\
\end{array}
\right)
\bigoplus
\left(
\begin{array}{cc}
-2 & 1\\
1 & -2\\
\end{array}
\right),
\end{equation}
hence $\disc(V_A)=-3$. It follows that the map
$V_A\ra\NS(A)$ is injective, and  we may regard $V_A$
as a subgroup of both $\Pic(A)$ or $\NS(A)$.

The exact sequence $0\ra\alpha_2\ra A\ra B\ra 0$ 
induces an exact sequence
$$
0\lra \Hom(\alpha_2,\GG_m)\lra \Pic^0_{B/k}\lra\Pic^0_{A/k}\lra 0
$$
by \cite{Mumford 1970}, Theorem 1 on page 143.
So $\Pic(B)\ra\Pic(A)$ is injective, 
and we can form the intersection $V_B=V_A\cap\Pic(B)$ inside $\Pic(A)$.
Note that we may regard $V_B$ as either a subgroup
$V_B\subset\Pic(B)$ or as a subgroup $V_B\subset\Pic(A)$.

\begin{proposition}
\mylabel{annihilated two}
The group $V_A/V_B$ is annihilated by 2.
\end{proposition}

\proof
Given an invertible  $\O_A$-module $\shL$, let $K(\shL)\subset A$ 
be the kernel of the morphism 
$A\ra\Pic^0_{A/k}$, $x\mapsto T_x^*(\shL)\otimes\shL^\vee$.
According to \cite{Moret-Bailly 1981},
Corollary 4.1.1 and Theorem 2.3 we have $\shL=g^*(\shN)$ 
for some invertible $\O_B$-module
$\shN$ if and only if $\alpha_2\subset K(\shL)$. For $\shL_i=\O_A(C_i)$ 
we have $K(\shL_i)=C_i$ and $K(\shL_i^{\otimes 2})=2C_i$,
by considering the quotient morphism $A\ra A/C_i$.
Obviously we have $\alpha_2\subset 2C_i$, so
$\shL_i^{\otimes 2}$ lies in $V_B$.
\qed

\medskip
The group $V_B$ is closely related to $\Pic(S)$ and $\Pic(B/G)$:

\begin{proposition}
\mylabel{factors}
The map $\Pic(B/G)\ra\NS(A)$ factors through $V_B\subset\NS(A)$.
\end{proposition}

\proof
We have $\rho(B/G)\leq 4$ because $\rho(S)\leq 22$ and 
there are $18$ exceptional curves for the resolution of
singularities $S\ra B/G$.
The elliptic curves $C_i\subset A$ are $G$-invariant, 
and $2C_i$ are $\alpha_2$-invariant.
Let $h:A\ra B/G$ be the canonical morphism and $D_i=h(C_i)$.
The $D_i\subset B/G$ are integral Weil divisors. Since the singularities on 
$B/G$ are rational double points of type $A_2$, the Weil
divisors  $3D_i$ are Cartier.
By construction, $h^*(3D_i)=n_iC_i$ for some integers  $n_i\geq 1$.
It follows that $\Pic(B/G)\ra\NS(A)$ factors over 
$\NS(A)\cap (V_A\otimes\QQ)$. 
Since $\disc(V_A)=-3$ is square free,
the subgroup $V_A\subset\NS(A)$ is a direct summand, so the map factors
over $V_A$,
and therefore over $V_B$ as well.
\qed

\medskip
The subgroup $V_B\subset V_A$ depends on the position of $\alpha_2\subset A$.
 Let $x\in A$ be the origin, and choose a regular system
of parameters $u,v\in\O_{A,x}^\wedge$ 
such that $u\mapsto u+u'$, $v\mapsto v+v'$ is the formal 
group law for the formal group $A^\wedge$.
Write $\alpha_2=\Spec k[\epsilon]$. 
Then the embedding 
$\alpha_2\ra A^\wedge$ is of the form $u\mapsto u+i\epsilon$ and 
$v\mapsto v+j\epsilon$ for 
some $(i,j)\in \AA^2(k)-\left\{0\right\}$. 
The subgroup scheme $\alpha_2\subset A$ 
depends only on the homogeneous coordinates
$(i:j)\in\PP^1(k)$. Let us regard the $\FF_4$-valued points
$\PP^1(\FF_4)$ as  a subset of the $k$-valued points $\PP^1(k)$.

\begin{proposition}
\mylabel{rank}
Set $v=\dim_{\FF_2}(V_A/V_B)$. Then $v=3$ if the homogeneous
coordinates satisfy $(i:j)\in\PP^1(\FF_4)$, 
and $v=4$ otherwise.
\end{proposition}

\proof
Consider an invertible $\O_A$-module of the form
$\shL=\O_A(\sum_{i=1}^4 n_iC_i')$,
where $n_i\in\left\{0,1\right\}$ and at least 
one coefficient $n_i$  is nonzero. 
The idea is to check 
whether $\alpha_2\subset A$ might be contained in  $K(\shL)\subset A$.
First note that $K(\shL)$ is a finite group scheme of 
length $\chi(\shL)^2$ 
if $\chi(\shL)\neq 0$,
and of dimension $>0$ otherwise,
according to \cite{Mumford 1970}, page 150.
The Riemann--Roch Theorem  and $n_i^2=n_i$ gives
$$
\chi(\shL)=(\sum_{i=1}^4 n_iC'_i)^2/2=n_1n_2-n_3-n_4+n_3n_4.
$$
If $\chi(L)$ is odd then $K(\shL)$ is an finite group scheme of odd length,
hence $\alpha_2\not\subset K(\shL)$.
The number $\chi(\shL)$ is odd except  if $(n_1,\ldots, n_4)$ is of the form
$$
(1,0,0,0),\quad (0,1,0,0),\quad (1,1,1,0),\quad (1,1,0,1),\quad (1,1,1,1).
$$
In these cases we compute the group schemes $K(\shL)\subset A$.
In the first four cases we have $\sum_{i=1}^4 n_iC_i'= C_j$ for some 
$1\leq j\leq 4$, and hence $K(\shL)=C_j$.

It remains to treat the last case $\sum_{i=1}^4C'_i$.
Modulo $2V_A\subset\Pic(A)$, 
this divisor is linearly equivalent to 
$$
C=C_1'+C_2'-3C_3'-3C_4',
$$
so $\sum_{i=1}^4C'_i\in V_B$ if and only if $C\in V_B$.
Consider the graph $C_5=\Gamma_{\varphi^2}$ of the automorphism 
$\varphi^2=\varphi^{-1}:E\ra E$.
Since $C_5$ is $G$-invariant and 
$2C_5$ is $\alpha_2$-invariant,
some multiple of $\O_A(C_5)\in\Pic(A)$ 
is the preimage of some class in $\Pic(B/G)$.
By Proposition \ref{factors}, we must have $C_5\in V_A\subset\NS(A)$.
But then $C$ must be numerically equivalent to $C_5$, because
$C\cdot C_i=C_5\cdot C_i$ for $1\leq i\leq 4$.
In turn, we have $K(\O_A(C))=K(\O_A(C_5))=C_5$.

Now we are done: Recall that $\varphi:E\ra E$ was defined on affine coordinates
by $(x,y)\mapsto(\zeta x,y)$. An easy calculation shows 
that  $(i:j)\in\PP^1(\FF_4)$ holds if and only if
$\alpha_2\subset A$ is contained in precisely one $C_j$ from
the elliptic curves
$C_1=E\times 0$, $C_2=0\times E$, $C_3=\Delta_E$, $C_4=\Gamma_\varphi$, 
and $C_5=\Gamma_{\varphi^2}$.
In this situation we have $V_B=2V_A+\ZZ C_j$ and  $v=3$.
If $(i:j)\not\in\PP^1(\FF_4)$, then $V_B=2V_A$ and  $v=4$.
\qed

\medskip
We come to the main result of this section:

\begin{theorem}
\mylabel{discriminant}
The group $\Pic(S)$ of the K3-surface $S$ is free of rank $\rho(S)=22$.
The discriminant $d\in\ZZ$ of the intersection pairing on $\Pic(S)$
is $d=-2^2$ if the homogeneous coordinates satisfy $(i:j)\in\PP^1(\FF_4)$, 
and $d=-2^4$ otherwise.
\end{theorem}

\proof
We have $d=-2^n$ 
for some integer $n\geq 0$ as explained in \cite{Artin 1974}, page 556.
To determine this integer it suffices to compute the discriminant
of the 2-adic intersection form $\Pic(S)\otimes\ZZ_2$ up to units.
The intersection matrix with respect to  the exceptional divisors for the
resolution of singularities $S\ra B/G$ is 
$$
\bigoplus_{i=1}^9 \left(\begin{array}{cc} -2&1\\1&-2\end{array}\right).
$$
Its discriminant is the unit $27\in\ZZ_2$, hence
$\Pic(B/G)\otimes\ZZ_2\subset\Pic(S)\otimes\ZZ_2$ 
is an orthogonal direct summand, and it suffices to compute the
discriminant of $\Pic(B/G)\otimes\ZZ_2$.

I claim that the inclusion $\Pic(B/G)\otimes\ZZ_2\subset V_B\otimes\ZZ_2$ 
from Proposition \ref{factors}
is bijective. 
Indeed: If $g:B\ra B/G$ denotes the projection map, we have
$g_*g^*(D)=3D$ for any  divisor $D\subset B/G$ by the projection formula.
This implies that the inclusion $\Pic(B/G)\subset V_B$, 
which has finite index,
becomes bijective after tensoring with $\ZZ_2$.
Moreover, the discriminants of the intersection forms on $\Pic(B/G)$
and $V_B\subset\Pic(B)$ differ only by a power of the 2-adic unit 
$3\in\ZZ_2$.

For our purpose it suffices to compute the discriminant of 
$V_B$ endowed with the intersection form from $V_B\subset\Pic(B)$.
Let us write $V_B'$ for the same module $V_B$, but endowed with
the intersection form from $V_B\subset\Pic(A)$.
Then $\disc(V_B)= 2^{-4}\disc(V_B')$. 
Indeed, if $h:A\ra B$ denotes the projection map, 
then $2(\shL\cdot\shL)=(\shL_A\cdot\shL_A)$ for
any $\shL\in\Pic(B)$.

This does it: We  have $\disc(V_A)=-3$
and $V_A/V_B\simeq(V_A\otimes\ZZ_2)/(V_B\otimes\ZZ_2)$.
According to \cite{Serre 1979}, Corollary on page 49, this gives
$$
\disc(V'_B\otimes\ZZ_2)=
\disc(V_A\otimes\ZZ_2)\cdot \ord(V_A/V_B)^2=-3\cdot 2^{2v},
$$
where $v=\dim_{\FF_2}(V_A/V_B)$.
Using Proposition \ref{rank}, we infer that up to 2-adic units 
$\Pic(B/G)$ has discriminant $2^{2\cdot3-4}=2^2$ if $(i:j)\in\PP^1(\FF_4)$,
and $2^{2\cdot 4-4}=2^4$ otherwise.
\qed

\medskip
For a  K3-surface in characteristic $p>0$ 
with Picard number $\rho=22$ and discriminant
$d=-p^{2\sigma_0}$, 
the integer $\sigma_0$ is called the \emph{Artin invariant}.
We see that our K3-surface $S$ has Artin invariant
$\sigma_0(S)=1$ or $\sigma_0(S)=2$.

\section{Nonrigidity}
\mylabel{nonrigidity}

Let $k$ be an algebraically closed field of characteristic 
$p=2$. Up to isomorphism, there is only one supersingular elliptic
curve $E$, which is given by the Weierstrass equation $x^3=y^2+y$.
Set $A=E\times E$ and $X'=A\times\PP^1$, and choose an exact sequence
$$
0\lra\O_{\PP^1}(-2)\stackrel{r,s}{\lra}
\O_{\PP^1}^{\oplus 2}\lra\O_{\PP^1}(2)\lra 0
$$
given by two homogeneous quadratic polynomials
$r,s\in H^0(\PP^1,\O_{\PP^1}(2))$ 
without common zeros.
As explained in Section \ref{abelian pencil}, this defines
a relative subgroup scheme $H\subset X'$, whose quotient $X=X'/H$ 
has $\omega_X=\O_X$.
As in the preceding section, the group $G=\ZZ/3\ZZ$ acts on $E$,
hence diagonally on $A$, and the induced fiberwise action on $X'\ra\PP^1$
descends to a fiberwise action on $X\ra\PP^1$. Set $Z=X/G$ and let
$Y\ra Z$ be the minimal resolution of singularities.
The same arguments as in Section \ref{kummer pencil} show that 
$\omega_Y=\O_Y$, and that the fibers $Y_t$ for the induced fibration
$f:Y\ra\PP^1$ are K3-surfaces with Picard number $\rho(Y_t)=22$.
So the smooth threefold $Y$ satisfies the assumptions of Section
\ref{calculation}, hence is a Calabi--Yau threefold in characteristic
$p=2$
without any formal lifting to characteristic zero.
In particular, the deformations over the $W(k)$-algebras are obstructed.
However, we shall see that there are many unobstructed deformations
over $k$-algebras.

\begin{proposition}
\mylabel{count invariant}
The fiber $Y_t$ has Artin invariant $\sigma_0(Y_t)=1$ if and only if
the homogeneous coordinates satisfy 
$(r(t):s(t))\in\PP^1(\FF_4)$, and $\sigma_0(Y_t)=2$ otherwise.
\end{proposition}

\proof
This follows from Proposition \ref{discriminant} and the definition
of $Y_y$ as the minimal resolution of the quotient of $X_t=X'_t/H_t$ 
by $G$.
\qed

\medskip
Choose homogeneous coordinates
$t=(t_0: t_1)$ for $\PP^1$ and consider, for example,
the case $r=t_0^2$ and $s=t_1^2$. Then
$(r:s)=(t_0^2:t_1^2)$ 
is $\FF_4$-valued if and only if $(t_0:t_1)$ is $\FF_4$-valued, and
there are precisely 5 such points on $\PP^1$.

Now consider the case that $r=t_0^2+at_0t_1$ and $s=t_1^2$, 
where $a\in k$ 
is a nonzero  element.
The quotient  $r/s=(t_0/t_1)^2+a(t_0/t_1)=P(t_0/t_1)$ 
is given by the separable polynomial $P(x)=x^2+ax$, hence
$P(x)=b$ has two distinct roots for each $b\in\FF_4$. 
This means that 
there are precisely 8 points on $\PP^1$ with $t_1\neq 0$
for which $(r(t):s(t))$ is 
$\FF_4$-valued, plus the additional point $(1:0)$.
 
This observation leads to nontrivial deformations over
the formal power series ring $R=k[[T]]$: Consider the two quadrics
$r=t_0^2+Tt_0t_1$ and $s=t_1^2$ with coefficients from $R$. 
It is easy to see that the construction of our
Calabi--Yau threefolds carries over to  a relative situation,
so the two quadrics $r,s$ define a family $Y\ra\Spec(R)$ of threefolds.
The second main result of this paper is that this are nontrivial
deformations:

\begin{theorem}
\mylabel{nontrivial}
The deformation $Y\ra\Spec(R)$ is not formally isomorphic to
the trivial deformation $Y_0\otimes R\ra\Spec(R)$.
\end{theorem}

\proof
Suppose to the contrary that these deformations are
formally isomorphic. Then they are isomorphic, and in particular
their generic geometric fibers $Y_L$ 
and $Y_0\otimes L$ are isomorphic, where $L$ is the algebraic
closure of $k((T))$.
But we saw above that the $Y_0\otimes L$ has 5  fibers over $\PP^1_L$
whose Artin invariant is $\sigma_0=1$, whereas $Y_L$ has 9 such fibers,
contradiction.
\qed

\section{The Artin--Mazur formal group}
\mylabel{artin mazur}

We now go back to  the situation in Section \ref{calculation},
such that $Y$ is a smooth proper threefold with $\omega_Y=\O_Y$ endowed with
a smooth morphism $f:Y\ra\PP^1$ whose geometric 
fibers are K3-surfaces with $\rho=22$
in characteristic $p>0$.
The goal of this section is to compute the \emph{Artin--Mazur formal
group} $\Phi^3_Y$ of such  Calabi--Yau threefolds.
Let me recall some definitions: If $A$ is a local Artin $k$-algebra
with maximal ideal $I\subset A$, we may view the sheaf 
$1+I\otimes_k\O_Y$
as an abelian subsheaf of $(A\otimes_k\O_Y)^\times$. This gives
a functor of Artin rings $\Phi_Y^3(A)=H^3(Y,1+I\otimes_k\O_Y)$.
Artin and Mazur proved in \cite{Artin; Mazur 1977}, Proposition 1.8
that, for schemes like $Y$, this functor
is representable by a 1-dimensional formal group.

\begin{proposition}
\mylabel{unipotent}
The Artin--Mazur formal group
$\Phi^3_Y$ of our Calabi--Yau threefold $Y$
is isomorphic to the formal additive group  $\hat{\GG}_a$.
\end{proposition}

\proof
The referee suggested the following direct argument:
Let $W=W(k)$ be the ring of Witt vectors, and $W\subset K$ 
be its quotient field.
It follows from \cite{Katz; Messing 1974}, Corollary 1 that the $l$-adic
Betti numbers $b_i$ for the projective smooth scheme $Y$ coincide
with the $W$-rank of the crystalline cohomology groups 
$H^i(Y/W,\O_{Y/W})$. In particular, $H^3(Y/W,\O_{Y/W})\otimes K=0$.
On the other hand, the part of the $F$-crystal 
$H^3(Y/W,\O_{Y/W})\otimes K$ whose
slopes fall into the interval $[0,1)$ 
is the Dieudonn\'e module for the maximal
$p$-divisible quotient of the formal group $\Phi^3_Y$
(compare \cite{Illusie 1979}, Remark 5.13 and
\cite{Artin; Mazur 1977}, Corollary 3.3).
It then follows from the classification of 1-dimensional formal groups
that $\Phi^3_Y=\hat{\GG}_a$.
\qed

\medskip
It follows from a result of Geer and Katsura (\cite{Geer; Katsura 2000},
Proposition 16.4)
that the Frobenius  operator $F$ on the Witt vector cohomology
$H^3(Y,W_i\O_Y)$ vanishes for all $i\geq 1$. Indeed, 
the Witt vector cohomology $H^3(Y,W\O_Y)$ 
is the Dieudonn\'e module for the formal group $\Phi^3_Y$.


\end{document}